\documentclass[11pt]{amsart}
\usepackage{ifthen,etoolbox,amsmath,amssymb,color,cite,mathtools,verbatim}
\usepackage[colorlinks,pagebackref,pdftex, bookmarks=false]{hyperref}
\apptocmd{\sloppy}{\hbadness 10000\relax}{}{}

\numberwithin{equation}{section}

\newtheorem{thm}[equation]{Theorem}
\newtheorem{prop}[equation]{Proposition}
\newtheorem{lemma}[equation]{Lemma}
\newtheorem{cor}[equation]{Corollary}

\theoremstyle{definition}
\newtheorem{rmk}[equation]{Remark}

\newtheorem{defn}[equation]{Definition}

\newcommand{\F}{\mathbb{F}}
\newcommand{\bP}{\mathbb{P}}
\newcommand{\R}{\mathcal R}
\newcommand{\Z}{\mathbb{Z}}
\newcommand{\Q}{\mathbb{Q}}
\newcommand{\C}{\mathbb{C}}

\DeclareMathOperator{\charp}{char}
\DeclareMathOperator{\PSL}{PSL}
\DeclareMathOperator{\PGammaL}{P\Gamma L}
\DeclareMathOperator{\AGL}{AGL}
\newcommand{\mybar}[1]{#1\llap{$\overline{\phantom{\rm#1}}$}}
\newcommand{\abs}[1]{\lvert #1 \rvert}

\begin{document}

\title{Determination of hyperovals by lines through a few points}

\author{Zhiguo Ding}
\address{
  Department of Mathematics,
  University of Michigan,
  530 Church Street,
  Ann Arbor, MI 48109-1043 USA
}
\email{dingz@umich.edu}

\author{Michael E. Zieve}
\address{
  Department of Mathematics,
  University of Michigan,
  530 Church Street,
  Ann Arbor, MI 48109-1043 USA
}
\email{zieve@umich.edu}
\urladdr{https://dept.math.lsa.umich.edu/$\sim$zieve/}


\date{\today}

\thanks{The authors thank Kai-Uwe Schmidt for providing valuable comments about this paper shortly before his untimely death.
The first author was supported in part by the Natural Science Foundation of Hunan Province of China (No.\ 2020JJ4164).}

\begin{abstract}
If $\mathcal{S}$ is a set of $q+2$ points in $\bP^2(\F_q)$ such that some point of $\mathcal{S}$ is not on any line containing two other points of $\mathcal{S}$, then in suitable coordinates $\mathcal{S}$ has the form $\mathcal{S}_f:=\{(c:f(c):1):c\in\F_q\}\cup\{(1:0:0), (0:1:0)\}$ for some $f(X)\in\F_q[X]$.
Let $\mathcal{T}$ be a subset of $\mathcal{S}_f$ which contains the two infinite points and at least $3+\log_3(q)/4$ finite points. We show that if there is no line passing through a point of $\mathcal{T}$ and two other points of $\mathcal{S}_f$, and $\deg(f)\le q^{1/4}$, then no three points of $\mathcal{S}_f$ are collinear, so that $\mathcal{S}_f$ is a hyperoval.
We also determine all $f(X)$ with $\deg(f)\le q^{1/4}$ for which $\mathcal{S}_f$ is a hyperoval, which strengthens a result that was proved by Caullery and Schmidt using entirely different methods.
\end{abstract}

\dedicatory{Dedicated to the memory of Kai-Uwe Schmidt}

\maketitle


\section{Introduction}

For any prime power $q$, a \emph{hyperoval} in $\bP^2(\F_q)$ is a set of $q+2$ points in $\bP^2(\F_q)$ which contains no three collinear points. It is well-known that there are no hyperovals when $q$ is odd \cite{Bose}. There has been a great deal of work on hyperovals in the past $75$ years, and several interesting examples have been found; see for instance the surveys \cite{Cherowitzo, Hirschfeld, Penttila}.

In this paper we show that, under certain hypotheses, knowing that any of a small set of lines contains no three points from a prescribed $(q+2)$-element set implies that no line whatsoever contains three such points, so that the set is a hyperoval.
Our results use the following notation.
\begin{defn}
For any $f(X)\in\F_q[X]$ we write $\mathcal{S}_{f,q}$ for the subset of $\bP^2(\F_q)$ defined by
\[
\mathcal{S}_{f,q}:=\{(c:f(c):1):c\in\F_q\}\cup\{(1:0:0),(0:1:0)\}.
\]
When the choice of $q$ is clear, we write $\mathcal{S}_{f,q}$ as $\mathcal{S}_f$.  
\end{defn}
We show in Lemma~\ref{o1} that, if $\mathcal S$ is a $(q+2)$-element subset $\bP^2(\F_q)$, and $\mathcal{S}$ contains a point $P$ which is not on any line containing two other points of $\mathcal{S}$, then there are coordinates for $\bP^2(\F_q)$ in which $\mathcal S$ has the form $\mathcal S_{f,q}$ for some $f(X)\in\F_q[X]$. Thus it suffices to analyze the sets $\mathcal{S}_{f,q}$. Our main result is as follows.

\begin{thm}\label{main}
Let $q=2^k$ with $k\ge 1$, and let $f(X)\in\F_q[X]$ have degree $n$, where $1<n\le q^{1/4}$. Write $N:=3+\log_3(n)$, and let $\mathcal{T}$ be a subset of $\mathcal{S}_f$ which contains $(1:0:0)$, $(0:1:0)$, and at least $N$ other points. Then the following are equivalent:
\begin{enumerate}
\item there is no line which contains a point of $\mathcal{T}$ and at least two other points of $\mathcal{S}_f$;
\item $\mathcal{S}_f$ is a hyperoval;
\item $f(X)=\rho(X)\circ X^n\circ\eta(X)$ for some degree-one $\rho(X),\eta(X)\in\F_q[X]$ where one of the following holds:
\begin{itemize}
\item $n=6$ and $k$ is an odd positive integer,
\item $n=2^{\ell}$ for some positive integer $\ell$ with $\gcd(k,\ell)=1$.
\end{itemize}
\end{enumerate}
\end{thm}

The most novel part of Theorem~\ref{main} is the implication that (1) implies (2). This says that if no line through one of a small number of points of $\mathcal{S}_f$ contains three points of $\mathcal{S}_f$, then no line whatsoever contains three points of $\mathcal{S}_f$. 
This is a new type of implication in this topic.
In Section~\ref{sec:heuristic} we present a heuristic suggesting that the bound $3+\log_3(n)$ in Theorem~\ref{main} would be roughly the best possible if we remove the condition $n\le q^{1/4}$.

Segre showed that (3) of Theorem~\ref{main} implies (2), even without assuming $n\le q^{1/4}$ \cite{Segre-arcs,Segre-char2}. However, there are several known examples showing that (2) would not imply (3) if we did not assume $n\le q^{1/4}$ \cite{Penttila}.

Theorem~\ref{main} refines \cite[Thm.~1.2]{CS}, which shows that (2) implies (3) when our hypothesis $n\le q^{1/4}$ is replaced by the slightly more restrictive hypothesis $n<\frac12q^{1/4}$. Our methods are completely different from those in \cite{CS}, so that in addition to the new implication that (1) implies (2), we also obtain a new proof of \cite[Thm.~1.2]{CS}.
One consequence of the fact that (2) implies (3) is that the functions in (3) are the only $f(X)\in\F_q[X]$ for which there are infinitely many $m$ such that $\mathcal{S}_{f,q^m}$ is a hyperoval in $\bP^2(\F_{q^m})$. The special case of this last implication in which $f(X)$ is a monomial was conjectured by Segre and Bartocci \cite{SB} and proved in \cite{HM,Z-planar}.

Our proof of Theorem~\ref{main} combines many different types of ingredients, including among other things functional decomposition of odd polynomials and a new result of independent interest (Corollary~\ref{EPbp}) bounding the number of critical values of a low-degree permutation polynomial.
Since our proof of the latter result ultimately relies on the classification of finite simple groups, we will sketch a variant of our argument which only uses elementary results, and which yields a variant of Theorem~\ref{main} with $N$ replaced by $n/2$. We note that this elementary variant still refines \cite[Thm.~1.2]{CS}, and also shows that (1) implies (2) when $N$ is replaced by $n/2$.

It is natural to seek an analogue of Theorem~\ref{main} in case $q$ is odd.  Although it is known that in this case there are no hyperovals, the following stronger result seems to be new.

\begin{prop}\label{oddq}
If $q$ is an odd prime power, and $\mathcal{S}$ is a $(q+2)$-element subset of\/ $\bP^2(\F_q)$, then there are at most two points $P$ in $\mathcal{S}$ with the property that no line through $P$ contains more than two points of $\mathcal{S}$.
\end{prop}

Our proof of Proposition~\ref{oddq} is short and self-contained, and uses the same strategy as the classical proof that a finite group with nontrivial cyclic Sylow $2$-subgroups has no complete mapping \cite[Thm.~5]{HP}.

This paper is organized as follows. In the next section we include various facts which are used in our proof of Theorem~\ref{main} and which follow relatively easily from known results or known arguments. Then in Section~\ref{sec:main} we prove Theorem~\ref{main}, and in Section~\ref{sec:elem} we sketch an elementary proof of a variant of Theorem~\ref{main}.
In Section~\ref{sec:heuristic} we present heuristics predicting that the value $N$ in Theorem~\ref{main} is essentially the best possible.  Finally, we prove Proposition~\ref{oddq} in Section~\ref{sec:odd}.


\section{Background material}

In this section we present various facts which follow easily from known facts or known arguments.

Throughout this paper, $q$ denotes a fixed prime power, and for any $f(X)\in\F_q[X]$ we write $\mathcal S_f$ for the subset of $\bP^2(\F_q)$ defined by
\[
\mathcal{S}_f := \{(c:f(c):1):c\in\F_q\}\cup\{(1:0:0),(0:1:0)\}.
\]
Also, if $K$ is a field then $\mybar K$ denotes the algebraic closure of $K$.


\subsection{Collinearity and permutation polynomials}

\begin{lemma}\label{o1}
If $\mathcal{S}$ is a subset of\/ $\bP^2(\F_q)$ with $\abs{\mathcal{S}}=q+2$, and there is some $P\in \mathcal{S}$ for which no line through $P$ contains at least three points of $\mathcal{S}$, then in suitable coordinates $\mathcal{S}$ is $\mathcal S_f$ for some $f(X)\in\F_q[X]$.  If there are two such points $P\in\mathcal S$ then $\mathcal S$ is $\mathcal S_f$ for some $f(X)\in\F_q[X]$ which permutes\/ $\F_q$.
\end{lemma}

\begin{proof}
Pick $Q\in\mathcal{S}\setminus\{P\}$, where if possible we choose $Q$ so that no line through $Q$ contains at least three points of $\mathcal S$.
Choose coordinates so that $P$ and $Q$ become $(0:1:0)$ and $(1:0:0)$, respectively. By hypothesis, $P$ and $Q$ are the only points of $\mathcal S$ lying on the line $Z=0$.
For $c\in\F_q$, the line $X=cZ$ contains $P$, so there is at most one $d\in\F_q$ for which $\mathcal{S}$ contains $(c:d:1)$. Since $\mathcal{S}\setminus\{P,Q\}$ has size $q$, and contains no points on $Z=0$, each $c\in\F_q$ must correspond to a unique $d\in\F_q$. Thus $\mathcal{S}$ has the form $\mathcal{S}_f$.  Finally, if no line through $Q$ contains at least three points of $\mathcal S$ then a similar argument shows that 
for each $d\in\F_q$ there is exactly one $c\in\F_q$ for which $\mathcal{S}$ contains $(c:d:1)$, so that
$f(X)$ permutes $\F_q$.
\end{proof}

\begin{lemma}\label{o2}
Suppose $q>2$ and $f(X)\in\F_q[X]$ satisfies $\deg(f)<q$. Then
\begin{enumerate}
\item each line through $(0:1:0)$ contains exactly two points of $\mathcal{S}_f$;
\item $f(X)$ permutes\/ $\F_q$ if and only if each line through $(1:0:0)$ contains at most two points of $\mathcal{S}_f$;
\item for each $a\in\F_q$, the polynomial $(f(X+a)-f(a))/X$ permutes\/ $\F_q$ if and only if each line which passes through $(a:f(a):1)$ contains at most two points of $\mathcal{S}_f\setminus\{(1:0:0)\}$.
\end{enumerate}
\end{lemma}

\begin{proof}
Item (1) holds because the lines through $(0:1:0)$ are the lines $Z=0$ and $X=dZ$ with $d\in\F_q$. Likewise (2) holds. For (3), note that for distinct $a,b\in\F_q$ the line through $(a:f(a):1)$ and $(b:f(b):1)$ is
\[
(b-a)\cdot\bigl(Y-f(a)Z\bigr) = \bigl(f(b)-f(a)\bigr)\cdot(X-aZ).
\]
For fixed $a$, this line is determined by the point $(b-a:f(b)-f(a))$ in $\bP^1(\F_q)$. Thus the line is uniquely determined by the choice of $b\in\F_q\setminus\{a\}$ if and only if the polynomial $(f(X)-f(a))/(X-a)$ is injective on $\F_q\setminus\{a\}$, or equivalently $g_a(X):=(f(X+a)-f(a))/X$ is injective on $\F_q^*$.
Since $\deg(g_a)<q-1$, we know that $\sum_{c\in\F_q} g_a(c)=0$ \cite[\S 9]{Dickson}. Since also $\sum_{c\in\F_q} c = 0$ (because $q>2$), therefore $g_a(X)$ is injective on $\F_q^*$ if and only if $g_a(X)$ permutes $\F_q$. This implies (3).
\end{proof}


\subsection{Decompositions of odd polynomials}

It is known that an odd polynomial in $\C[X]$ can only decompose as the composition of two odd polynomials, up to linear changes of variables (cf.\ \cite[Prop.~1]{HR} or \cite[Lemma~3.11]{ZM}). We will use a generalization of this result to polynomials over an arbitrary field $K$. 
This generalization can be proved by the same argument as the case $K=\C$, although some care (and a new type of hypothesis) is needed in this more general situation. For the reader's convenience, we include the proof of this generalization in order to clarify the new issue occurring in positive characteristic.

\begin{defn}
If $K$ is a field and $\rho(X)\in K[X]$ has degree one, then we write $\rho^{-1}(X)$ for the unique degree-one polynomial in $K[X]$ such that $\rho^{-1}(\rho(X))=X=\rho(\rho^{-1}(X))$. Explicitly, if $\rho(X)=aX+b$ then $\rho^{-1}(X):=(X-b)/a$.
\end{defn}

\begin{lemma}\label{odd}
Let $K$ be a field, and let $f(X)\in X K[X^2]$ be a monic polynomial. If $f(X)=g(X)\circ h(X)$ with $g,h\in K[X]$, where $\charp(K)\nmid\deg(g)$, then there exists a degree-one $\rho(X)\in K[X]$ such that both $g(\rho(X))$ and $\rho^{-1}(h(X))$ are monic polynomials in $X K[X^2]$.
\end{lemma}

\begin{proof}
Denote the leading term of $h(X)$ by $aX^n$, and write $b:=h(0)$. Define $\rho(X):=aX+b$, so that $\widetilde h(X):=\rho^{-1}(h(X))$ is a monic polynomial with no constant term. Then $\widetilde g(X):=g(\rho(X))$ satisfies $\widetilde g(X)\circ\widetilde h(X)=g(X)\circ h(X)=f(X)$. 
Since the leading term of $\widetilde h(X)$ is $X^n$, if the leading term of $\widetilde g(X)$ is $cX^m$ then the leading term of $f(X)$ is $cX^{mn}$, so since $f(X)$ is monic we must have $c=1$. Moreover, since $f(X)\in XK[X^2]$, in particular $\deg(f)=mn$ is odd, so that $m$ and $n$ are odd.

First suppose that $\widetilde h(X)$ has a term of even degree, and let $dX^\ell$ be the highest-degree term of $\widetilde h(X)$ having even degree. Since $\widetilde h(0)=0$, we know that $\ell>0$. Since $m$ is odd, the $m$-th power of a polynomial in $X K[X^2]$ will also be in $X K[X^2]$. 
Thus the highest-degree term of $h(X)^m$ having even degree is $mdX^{(m-1)n+\ell}$, where $md\ne 0$ due to the hypothesis $\charp(K)\nmid\deg(g)=m$. Since $f(X)-h(X)^m$ has degree at most $(m-1)n$, which is less than $(m-1)n+\ell$, we conclude that $mdX^{(m-1)n+\ell}$ is a term of $f(X)$, contradicting the hypothesis $f(X)\in X K[X^2]$. This contradiction implies that $\widetilde h(X)\in X K[X^2]$.

Next suppose that $\widetilde g(X)$ has a term of even degree, and let $eX^k$ be the highest-degree term of $\widetilde g(X)$ having even degree. Since the composition of two polynomials in $X K[X^2]$ will also be in $X K[X^2]$, it follows that $eX^{kn}$ is a term of $f(X)$, contradicting $f(X)\in X K[X^2]$.  Thus $\widetilde g(X)\in X K[X^2]$, which concludes the proof.
\end{proof}

\begin{rmk}
Lemma~\ref{odd} would not be true without the hypothesis on $\deg(g)$: for instance, if $K=\F_3$ then $(X^3-X^2+X)\circ (X^3+X^2+X)=X^9+X^5+X$.
\end{rmk}

\begin{cor}\label{oddcor}
Let $K$ be a field, and let $f(X)\in X K[X^2]$ be a monic polynomial with $\charp(K)\nmid\deg(f)$. If $f(X)=u_1(X)\circ u_2(X)\circ\dots\circ u_k(X)$ with $u_i(X)\in K[X]$, then there exist degree-one $\rho_1(X),\rho_2(X),\dots,\rho_{k-1}(X)\in K[X]$ such that, if we write $\rho_0(X)=\rho_k(X)=X$, then for $1\le i\le k$ the polynomial $\rho_{i-1}^{-1}(X)\circ u_i(X)\circ\rho_i(X)$ is monic and is in $X K[X^2]$.
\end{cor}

\begin{proof}
We proceed by induction on $k$, with the base case $k=1$ being vacuously true.  For the inductive step, write $g(X):=u_1(X)\circ\dots\circ u_{k-1}(X)$. By Lemma~\ref{odd}, there is a degree-one $\rho_{k-1}(X)\in K[X]$ such that both $\widetilde g(X):=g(\rho_{k-1}(X))$ and $\rho_{k-1}^{-1}(u_k(X))$ are monic polynomials in $X K[X^2]$. 
Now the result follows by induction, since $\widetilde g(X)=u_1(X)\circ\dots\circ u_{k-2}(X)\circ \widetilde u_{k-1}(X)$ where $\widetilde u_{k-1}(X):=u_{k-1}(\rho_{k-1}(X))$.
\end{proof}


\subsection{Separable polynomials}

We begin with the following immediate observation.

\begin{lemma}\label{sep}
For any field $K$ of characteristic $p$, and any $f(X)\in K[X]$, the following are equivalent:
\begin{enumerate}
\item the derivative $f'(X)$ is the zero polynomial;
\item $f(X)\in K[X^p]$.
\end{enumerate}
\end{lemma}

\begin{defn}
For any field $K$, we say that $f(X)\in K[X]$ is \emph{separable} if it does not satisfy the two equivalent conditions in Lemma~\ref{sep}.
\end{defn}

\begin{rmk}
We caution the reader that this definition differs from familiar definitions of separability of a polynomial. The motivation for this definition is that, in the case of nonconstant $f(X)$, it is equivalent to separability of the field extension $K(x)/K(f(x))$, where $x$ is transcendental over $K$. 
Likewise, if $f(X)$ is nonconstant then our definition is equivalent to asserting that $f(X)-t$ has no multiple roots in any extension of $K(t)$, where $t$ is transcendental over $K$. These equivalences are shown, for instance, in \cite[Lemma~2.2]{DZ4}.
\end{rmk}


\subsection{Critical points and critical values}

\begin{defn}
For any field $K$ and any nonconstant $f(X)\in K[X]$, the \emph{critical points} of $f(X)$ are the elements $c\in \mybar K$ for which $f'(c)=0$, and the \emph{critical values} of $f(X)$ are the values $f(c)$ where $c$ is a critical point of $f(X)$.
\end{defn}

\begin{lemma}\label{rh}
If $K$ is a field of characteristic $2$, then any separable $f(X)\in K[X]$ of degree $n>0$ has at most $(n-1)/2$ critical points.
\end{lemma}

\begin{proof}
The critical points are the roots of $f'(X)$, which is a polynomial of degree at most $n-1$ having no odd-degree terms, so that $f'(X)=g(X^2)$ for some $g(X)\in K[X]$ with $\deg(g)\le (n-1)/2$. 
Since $f(X)$ is separable, we have $f'(X)\ne 0$, so that $g(X)\ne 0$.  Thus $g(X)$ has at most $\deg(g)$ roots, and the square roots of these roots are the roots of $f'(X)$.
\end{proof}


\subsection{Ramification}

\begin{defn}
For any field $K$, any nonconstant $f(X)\in K[X]$, and any $c,d\in\mybar K$, the \emph{ramification index} $e_f(c)$ is the multiplicity of $c$ as a root of $f(X)-f(c)$, and the \emph{ramification multiset} $\R_f(d)$ is the multiset of multiplicities of roots of $f(X)-d$, or equivalently, the multiset of all $e_f(b)$ with $b$ varying over all $f$-preimages of $d$ in $\mybar K$.
\end{defn}

We will use the following basic facts.

\begin{lemma}\label{ram}
Let $K$ be a field, and pick nonconstant $f,g\in K[X]$ and any $b,c,d\in\mybar K$. Then
\begin{enumerate}
\item $\R_f(d)$ is a collection of positive integers whose sum is $\deg(f)$;
\item $c$ is a critical point of $f(X)$ if and only if $e_f(c)>1$;
\item $e_{f\circ g}(b)=e_f(g(b))\cdot e_g(b)$.
\end{enumerate}
\end{lemma}

\begin{proof}
Item (1) follows from unique prime factorization in $\mybar K[X]$. Item (2) is a restatement of the well-known fact that $c$ is a multiple root of $f(X)-f(c)$ if and only if $c$ is a root of the derivative of this polynomial. 
Finally, item (3) is well-known, and in this case can be seen as follows: $e_g(b)$ is the degree of the lowest-degree term of $\widetilde{g}(X):=g(X+b)-g(b)$, and $e_f(g(b))$ is the degree of the lowest-degree term of $\widetilde{f}(X):=f(X+g(b))-f(g(b))$, so that $e_f(g(b))\cdot e_g(b)$ is the degree of the lowest-degree term of
\[
\widetilde{f}(\widetilde{g}(X))=f(g(X+b))-f(g(b)),
\]
and hence equals $e_{f\circ g}(b)$.
\end{proof}

\begin{rmk}
When writing multisets, we use exponents to denote multiplicities, and we write an exponent of $*$ for a multiplicity which is not being written explicitly for typographical reasons. Thus, for instance, the assertion $\R_f(d)=[1,2^*]$ says that $n:=\deg(f)$ is odd and $\R_f(d)$ consists of a single $1$ and $(n-1)/2$ copies of $2$.
\end{rmk}


\subsection{Dickson polynomials}

\begin{defn}
For any field $K$, any element $a\in K$, and any positive integer $n$, the \emph{Dickson polynomial} of degree $n$ with parameter $a$ is the unique polynomial $D_n(X,a) \in K[X]$ satisfying the functional equation 
\begin{equation}\label{dickson}
D_n\Bigl(X+\frac{a}{X},a\Bigr) = X^n+\frac{a^n}{X^n};
\end{equation}
explicitly, we have
\begin{equation}\label{dicksoncoeffs}
D_n(X,a)=\sum_{i=0}^{\lfloor n/2\rfloor}\frac{n}{n-i}\binom{n-i}i (-a)^i X^{n-2i},
\end{equation}
where $\frac{n}{n-i}\binom{n-i}i$ is calculated in\/ $\Q$, and turns out to be in\/ $\Z$ (and hence can be interpreted as an element of $K$).
\end{defn}

\begin{rmk}
For proofs that the polynomial in \eqref{dicksoncoeffs} satisfies \eqref{dickson}, see \cite{ACZ,LMT}.
\end{rmk}

We record the following immediate consequence of this definition.

\begin{cor}\label{Dicksonscale}
For any field $K$, any $a,b\in K$, and any positive integer $n$, we have
\[
D_n(bX,b^2 a)=b^n D_n(X,a).
\]
\end{cor}

The following result is a slight generalization of \cite[\S 54]{Dickson}; see e.g.\ \cite[Thms.~3.1 and 3.2]{LMT}.

\begin{lemma}\label{DicksonPP}
Let $n$ be a positive integer, let $q$ be a prime power, and pick $a\in\F_q$. If $a=0$ then $D_n(X,a)=X^n$, which permutes\/ $\F_q$ if and only if $\gcd(n,q-1)=1$. If $a\ne 0$ then $D_n(X,a)$ permutes\/ $\F_q$ if and only if $\gcd(n,q^2-1)=1$. 
\end{lemma}

We now determine the ramification of Dickson polynomials over $\mybar\F_2$.

\begin{lemma}\label{Dicksonram}
Pick any $a\in\mybar\F_2^*$, and let $n$ be an odd integer with $n>1$. Then the following statements hold:
\begin{itemize}
\item $0$ is the unique critical value of $X^n$, and $\R_{X^n}(0)=[n]$ where $0$ is the unique root of $X^n$;
\item $0$ is the unique critical value of $f(X):=D_n(X,a)$, and $\R_f(0)=[1,2^*]$ where $0$ is the unique root of $f(X)$ with ramification index $1$.
\end{itemize}
\end{lemma}

\begin{proof}
The assertions about $X^n$ are immediate. Since $f(X):=D_n(X,a)$ satisfies
\[
f\Bigl(X+\frac{a}X\Bigr)=X^n+\frac{a^n}{X^n},
\]
taking derivatives yields
\[
\Bigl(1+\frac{a}{X^2}\Bigr)\cdot
f'\Bigl(X+\frac{a}X\Bigr) = X^{n-1}+\frac{a^n}{X^{n+1}}=X^{-n-1}(X^{2n}+a^n).
\]
Writing $b:=\sqrt{a}$, it follows that
\[
f'\Bigl(X+\frac{a}X\Bigr) = X^{-n-3}\cdot\frac{X^{2n}+a^n}{X^2+a} = nX^{-n-3}\cdot\prod_{\zeta\in\mu_n\setminus\{1\}} (X+b\zeta)^2,
\]
where we write $\mu_n$ for the set of $n$-th roots of unity in $\mybar\F_2$.
Thus the critical points of $f(X)$ are the elements $b(\zeta+\zeta^{-1})$ with $\zeta\in\mu_n\setminus\{1\}$.  In particular, $f(X)$ has exactly $(n-1)/2$ critical points.  For any $\zeta\in\mu_n$, we compute
\[
f\bigl(b(\zeta+\zeta^{-1})\bigr) = f\Bigl(b\zeta+\frac{a}{b\zeta}\Bigr) = b^n\zeta^n+\frac{a^n}{b^n\zeta^n},
\]
which equals $0$ since $\zeta^n=1$ and $a=b^2$. Thus $0$ is the unique critical value of $f(X)$, and also $0=b(1+1^{-1})$ is a root of $f(X)$ which is not a critical point of $f(X)$. Finally, Lemma~\ref{ram} implies that $\R_f(0)=[1,2^{(n-1)/2}]$.
\end{proof}

For more about Dickson polynomials, see \cite[p.~2872]{ACZ} or \cite{LMT}.


\subsection{Exceptional polynomials, I: elementary results}

\begin{defn}\label{EP}
A polynomial $f(X)\in\F_q[X]$ is \emph{exceptional over\/ $\F_q$} if the only polynomials in\/ $\F_q[X,Y]$ which divide $f(X)-f(Y)$ and are irreducible in\/ $\mybar\F_q[X,Y]$ are $c\cdot (X-Y)$ with $c\in\F_q^*$.
\end{defn}

The interesting aspect of exceptional polynomials is their connection with permutation polynomials, via the following two results.

\begin{lemma}\label{EPPP}
If $f(X)\in\F_q[X]$ is exceptional over\/ $\F_q$, then $f(X)$ permutes\/ $\F_q$. 
\end{lemma}

\begin{lemma}\label{PPEP}
If $f(X)\in\F_q[X]$ permutes\/ $\F_q$, and $\deg(f)\le q^{1/4}$, then $f(X)$ is exceptional over\/ $\F_q$.
\end{lemma}

\begin{rmk}
Lemma~\ref{EPPP} was first shown in \cite[Thm.~5]{Cohen1970}; see also \cite[Thm.~1.1]{GTZ}. The simplest correct proof we have seen for Lemma~\ref{PPEP} is in \cite{Fan}. A weaker version of Lemma~\ref{PPEP} was shown in \cite{BSD}, in which the inequality $\deg(f)\le q^{1/4}$ was replaced by the condition that $q$ is sufficiently large compared to $\deg(f)$.
Lemma~\ref{PPEP} can be proved by making explicit each step of the argument in \cite{BSD}, primarily by using the version of Weil's bound for singular affine plane curves in \cite[Cor.~2(b)]{LY}. We note that the proofs of Lemma~\ref{PPEP} in \cite[\S 3.2]{AP} and \cite[Thm.~1]{vzGvalues} are not correct, as explained in \cite[p.~697]{MullerZplanar} and \cite[Rmk.~8.4.20]{ZEP}.
\end{rmk}

We will use the following characterizations of exceptional polynomials.

\begin{prop}\label{PPEPequiv}
For $f(X)\in\F_q[X]$, the following are equivalent:
\begin{enumerate}
\item $f(X)$ is exceptional over\/ $\F_q$;
\item $f(X)$ permutes $\F_{q^m}$ for infinitely many positive integers $m$;
\item there is a positive integer $N$ such that $f(X)$ permutes\/ $\F_{q^m}$ for each positive integer $m$ with $\gcd(m,N)=1$.
\end{enumerate}
\end{prop}

The proof of Proposition~\ref{PPEPequiv} relies on the following consequence of Definition~\ref{EP}.

\begin{lemma}\label{dense}
If $f(X)\in\F_q[X]$ is exceptional over\/ $\F_q$, then there exists a positive integer $N$ such that $f(X)$ is exceptional over\/ $\F_{q^m}$ for every positive integer $m$ with $\gcd(m,N)=1$. 
\end{lemma}

\begin{proof}
Plainly exceptionality implies that $f(X)$ is nonconstant. Write $f(X)-f(Y) = a\cdot\prod_{i=1}^k H_i(X,Y)$ where $a\in\mybar\F_q^*$ and each $H_i(X,Y)$ is an irreducible polynomial in $\mybar\F_q[X,Y]$ in which at least one coefficient is in $\F_q^*$. 
Let $\Lambda$ be the set of all coefficients of all the polynomials $H_i(X,Y)$, so that $\Lambda$ is finite. Then the field $\F_q(\Lambda)$ equals $\F_{q^N}$ for some $N$. 
Pick $m$ with $\gcd(m,N)=1$, so that $\F_{q^m}\cap\F_{q^N}=\F_q$. Then any polynomial in $\F_{q^m}[X,Y]$ which divides $f(X)-f(Y)$ and is irreducible in $\mybar\F_q[X,Y]$ must be $b\cdot H_i(X,Y)$ for some $i$ and some $b\in\F_{q^m}^*$. 
Since the coefficients of $H_i(X,Y)$ are in both $\F_{q^m}$ and $\F_{q^N}$, they must be in $\F_q$, so that exceptionality of $f(X)$ over $\F_q$ implies $H_i(X,Y)=c\cdot (X-Y)$ with $c\in\F_q^*$. Thus $f(X)$ is exceptional over $\F_{q^m}$.
\end{proof}

\begin{rmk}
Lemma~\ref{dense} was first proved in \cite[Exceptionality Lemma, p.~167]{FGS}, via a different argument.
\end{rmk}

We also record the following immediate observation.

\begin{lemma}\label{triv}
If $f(X)\in\F_q[X]$ is exceptional over\/ $\F_{q^m}$ for some positive integer $m$, then $f(X)$ is exceptional over\/ $\F_q$.
\end{lemma}

\begin{proof}[Proof of Proposition~\ref{PPEPequiv}]
Plainly (3) implies (2). If (2) holds then Lemmas~\ref{PPEP} and \ref{triv} imply (1). Finally, if (1) holds then Lemmas~\ref{dense} and \ref{EPPP} yield (3).
\end{proof}

\begin{cor}\label{factor}
Pick $u_1(X),u_2(X),\dots,u_k(X)\in\F_q[X]$, and write $f(X):=u_1(X)\circ u_2(X)\circ\cdots\circ u_k(X)$. Then $f(X)$ is exceptional over\/ $\F_q$ if and only if each $u_i(X)$ is exceptional over\/ $\F_q$.  
\end{cor}

\begin{proof}
This follows from repeated application of Proposition~\ref{PPEPequiv}. If $f(X)$ is exceptional over $\F_q$, then $f(X)$ permutes $\F_{q^m}$ for infinitely many $m$, so that each $u_i(X)$ permutes each of these infinitely many $\F_{q^m}$ and hence is exceptional over $\F_q$.  
Conversely, if each $u_i(X)$ is exceptional over $\F_q$, then for each $i$ there is a positive integer $N_i$ such that $u_i(X)$ permutes $\F_{q^m}$ whenever $\gcd(m,N_i)=1$. Writing $N:=N_1N_2\cdots N_k$, it follows that $f(X)$ permutes $\F_{q^m}$ whenever $\gcd(m,N)=1$, so that $f(X)$ is exceptional over $\F_q$.
\end{proof}

\begin{rmk}
The ``if'' implication of Corollary~\ref{factor} was first shown in \cite[p.~60]{DL}. The full result was asserted in \cite[Thm.~1]{Fried-generalizedSchur}, but the proof there is incorrect, since it implicitly assumes that all polynomials $u_i(X)$ permute the same set of fields $\F_{q^m}$.  
The first correct proof of Corollary~\ref{factor} is \cite[Lemma~2.5]{Cohen-primitive}. See also \cite[Exceptionality Lemma, p.~167]{FGS} or \cite[Lemma~3.3(c)]{GW}. Our proof is simpler and more elementary than all previous correct proofs of either implication of Corollary~\ref{factor}.
\end{rmk}

We also need the following simple facts about low-degree exceptional polynomials.

\begin{lemma}\label{34}
Let $f(X)\in\F_q[X]$ be separable and exceptional. Then
\begin{itemize}
\item $\deg(f)\ne 2$;
\item if $\deg(f)=3$ then $f(X)$ has at most one critical value; 
\item if $\deg(f)=4$ then $q$ is even and $f(X)$ has no critical values.
\end{itemize}
\end{lemma}

Lemma~\ref{34} follows at once from the combination of Proposition~\ref{PPEPequiv} and the following classical result from \cite[\S 21--23 and 69--71]{Dickson}.

\begin{lemma}\label{dickson234}
Assume $q>7$ and $f(X)\in\F_q[X]$ has degree in $\{2,3,4\}$. Then $f(X)$ permutes\/ $\F_q$ if and only if there exist degree-one $\rho(X),\eta(X)\in\F_q[X]$ such that $\rho(X)\circ f(X)\circ\eta(X)$ is one of the following:
\begin{itemize}
\item $X^2$, if $q$ even;
\item $X^3$, if $q\not\equiv 1\pmod 3$;
\item $X^3-aX$, if $3\mid q$ and $a$ is a nonsquare in\/ $\F_q$;
\item $X^4+aX^2+bX$, if $q$ even and $X^3+aX+b$ has no roots in\/ $\F_q^*$.
\end{itemize}
\end{lemma}

\begin{defn}
For any field $K$, a polynomial $f(X)\in K[X]$ is \emph{indecomposable} if $\deg(f)>1$ and $f(X)$ cannot be written as the composition of lower-degree polynomials in $K[X]$.
\end{defn}

\begin{lemma}\label{tameEP}
If $f(X)\in\F_q[X]$ is indecomposable and exceptional, and $\deg(f)$ is coprime to $q$, then $f(X)=\rho(X)\circ g(X)\circ\eta(X)$ for some degree-one $\rho(X),\eta(X)\in\F_q[X]$ and some $g(X)$ which is either $X^n$ with $n$ an odd prime or $D_n(X,a)$ with $n\ge 5$ prime and $a\in\F_q^*$.
\end{lemma}

\begin{rmk}
The first correct proof of Lemma~\ref{tameEP} in the literature is  \cite[Appendix]{Muller-Schur} (in combination with the well-known Lemma~\ref{DicksonPP}). Earlier proofs in \cite[p.~177]{FGS} and \cite[Thm.~1]{Klyachko} were not correct, as explained in \cite[pp.~354--355]{TSchur}.
The proof relies on the Riemann--Hurwitz genus formula and classical group-theoretic results due to Burnside and Schur, for which short self-contained proofs appear in \cite{LMT}.
\end{rmk}

We will use the following combination of several results in this section.

\begin{cor}\label{oddEP}
Let $q$ be a prime power, and let $f(X)$ be a monic polynomial in $X\F_q[X^2]$ such that $\gcd(q,\deg(f))=1$ and $1<\deg(f)\le q^{1/4}$. If $f(X)$ permutes\/ $\F_q$ then $f(X)=u_1(X)\circ u_2(X)\circ\dots\circ u_k(X)$ where each $u_i(X)$ is either $X^{n_i}$ for some odd prime $n_i$, or $D_{n_i}(X,a_i)$ for some $a_i\in\F_q^*$ and some prime $n_i$ with $n_i\ge 5$.
\end{cor}

\begin{proof}
Write $f(X)=u_1(X)\circ u_2(X)\circ\dots\circ u_k(X)$ where each $u_i(X)$ is an indecomposable polynomial in $\F_q[X]$. By Corollary~\ref{oddcor} there is a choice of such $u_i(X)$'s for which each $u_i(X)$ is monic and is in $X\F_q[X^2]$. 
Since $f(X)$ permutes $\F_q$, also each $u_i(X)$ permutes $\F_q$, which by Lemma~\ref{PPEP} implies that each $u_i(X)$ is exceptional. By Lemma~\ref{tameEP} and Corollary~\ref{Dicksonscale}, each $u_i(X)$ is either $X^{n_i}$ with $n_i$ an odd prime or $D_{n_i}(X,a_i)$ with $n_i\ge 5$ prime and $a_i\in\F_q^*$.
\end{proof}


\subsection{Exceptional polynomials, II: hard results}

We now present some results about exceptional polynomials whose proofs depend on the classification of finite simple groups.

\begin{defn}
If $K$ is a field and $f(X)\in K[X]$ is separable, then the \emph{geometric monodromy group} of $f(X)$ is the Galois group of $f(X)-t$ over $\mybar K(t)$, where $t$ is transcendental over $K$.
\end{defn}

\begin{lemma}\label{FGS}
Let $q$ be a power of a prime $p$, and let $f(X)\in\F_q[X]$ be separable, indecomposable, and exceptional. Then $n:=\deg(f)$ and the geometric monodromy group $G$ of $f(X)$ satisfy one of the following:
\begin{enumerate}
\item $n$ is an odd prime with $n\ne p$ and $G\le\AGL_1(n)$;
\item $n=p^k$ and $G$ is a subgroup of $\AGL_k(p)$ such that the group of translations in $G$ acts transitively;
\item $p\in\{2,3\}$, $n=p^k(p^k-1)/2$ for some odd $k>1$, and $\PSL_2(p^k)\le G\le\PGammaL_2(p^k)$.
\end{enumerate}
\end{lemma}

\begin{proof}
This is the combination of \cite[Thm.~13.6 and Thm.~14.1]{FGS} with well-known results which can be found, for instance, in \cite[Lemma~2]{FriedSchur} and \cite[Exceptionality lemma]{FGS}.
\end{proof}

\begin{cor}\label{indEPbp}
If $q$ is an even prime power, and $f(X)\in\F_q[X]$ is separable, indecomposable, and exceptional, then $f(X)$ has at most one critical value.
\end{cor}

\begin{proof}
In case (1) of Lemma~\ref{FGS}, this follows from Lemmas~\ref{Dicksonram} and \ref{tameEP}. In cases (2) and (3) of Lemma~\ref{FGS}, it is shown in \cite[Lemma~2.1]{GM} and \cite[Thm.~2.1]{GZ}, respectively.
\end{proof}

\begin{cor}\label{EPbp}
If $q$ is an even prime power, and a degree-$n$ polynomial $f(X)\in\F_q[X]$ is separable and exceptional, then the number of critical values of $f(X)$ is at most $\log_3(n)$.
\end{cor}

\begin{proof}
Write $f(X)=\rho\circ u_1\circ u_2\circ\dots\circ u_k$ for some indecomposable polynomials $u_i(X)\in\F_q[X]$ and some $\rho(X)\in\F_q[X]$ of degree $1$. We prove the result by induction on $k$, with the base case $k=0$ being vacuous.
Write $g(X):=\rho\circ u_1\circ u_2\circ\dots\circ u_{k-1}$, so that $f(X)=g\circ u_k$, which implies the number of critical values of $f(X)$ is at most the sum of the numbers of critical values of $g(X)$ and $u_k(X)$. The hypothesis that $f(X)$ is separable implies that both $g(X)$ and $u_k(X)$ are separable. 
By Corollary~\ref{factor}, the hypothesis that $f(X)$ is exceptional implies that both $g(X)$ and $u_k(X)$ are exceptional. By the inductive hypothesis, $g(X)$ has at most $\log_3(\deg(g))$ critical values. Corollary~\ref{indEPbp} implies that $u_k(X)$ has at most one critical value.
By Lemma~\ref{34}, any separable indecomposable exceptional polynomial has degree at least $3$. Thus $\deg(g)\le n/3$, so the number of critical values of $g(X)$ is at most $\log_3(n/3)$, whence the number of critical values of $f(X)$ is at most $\log_3(n/3)+1=\log_3(n)$.
\end{proof}

%

\begin{rmk}
Minor variants of the above arguments show that if $q$ is odd and a degree-$n$ polynomial $f(X)\in\F_q[X]$ is separable and exceptional, then the number $N$ of critical values of $f(X)$ is at most $\log_{\sqrt{5}}(n)$, and also $N\le 2$ if in addition $f(X)$ is indecomposable.
\end{rmk}

\section{Proof of Theorem~\ref{main}}\label{sec:main}

In this section we prove Theorem~\ref{main}. In light of Lemma~\ref{o2}, it suffices to show the following result.

\begin{thm}\label{main2}
Let $q:=2^k$ with $k\ge 1$, and let $f(X)\in\F_q[X]$ have degree $n$, where $1<n\le q^{1/4}$. Write $N:=3+\log_3(n)$, and for any $a\in\F_q$ write $g_a(X):=\displaystyle{\frac{f(X+a)+f(a)}X}$. Then the following are equivalent:
\begin{enumerate}
\item $f(X)$ permutes\/ $\F_q$ and that there are at least $N$ elements $a\in\F_q$ for which $g_a(X)$ permutes\/ $\F_q$;
\item $f(X)$ permutes\/ $\F_q$ and $g_a(X)$ permutes\/ $\F_q$ for each $a\in\F_q$;
\item $f(X)=\rho(X)\circ X^n\circ\eta(X)$ for some degree-one $\rho(X),\eta(X)\in\F_q[X]$ where one of the following holds:
\begin{itemize}
\item $n=6$ and $k$ is odd,
\item $n=2^{\ell}$ for some positive integer $\ell$ with $\gcd(k,\ell)=1$.
\end{itemize}
\end{enumerate}
\end{thm}

We begin with the following known result, which is essentially due to Segre. We include a short proof for the reader's convenience.

\begin{lemma}\label{easy}
Let $q=2^k$ for some positive integer $k$. Then 
\begin{itemize}
\item $X^6$ permutes\/ $\F_q$ if and only if $k$ is odd, which holds if and only if $((X+a)^6+a^6)/X$ permutes\/ $\F_q$ for every $a\in\F_q$;
\item if $n=2^\ell$ for some positive integer $\ell$ then $X^n$ permutes\/ $\F_q$, and for each $a\in\F_q$ the polynomial $((X+a)^n+a^n)/X$ permutes\/ $\F_q$ if and only if $\gcd(k,\ell)=1$.
\end{itemize}
\end{lemma}

\begin{proof}
For $a\in\F_q$, the polynomial $g(X):=((X+a)^6+a^6)/X$ equals $D_5(X,a^2)$, so Lemma~\ref{DicksonPP} implies that if $a\ne 0$ then $g(X)$ permutes $\F_q$ if and only if $k$ is odd. Moreover, $X^6$ permutes $\F_q$ if and only if $k$ is odd, and if $k$ is odd then $D_5(X,0)=X^5$ permutes $\F_q$. 
Now write $n:=2^\ell$ with $\ell>0$, so that plainly $X^n$ permutes $\F_q$. For $a\in\F_q$, the polynomial $h(X):=((X+a)^n+a^n)/X$ equals $X^{n-1}$, which permutes $\F_q$ if and only if $\gcd(q-1,n-1)=1$; since $\gcd(q-1,n-1)=2^{\gcd(k,\ell)}-1$, it follows that $h(X)$ permutes $\F_q$ if and only if $\gcd(k,\ell)=1$.
\end{proof}

\begin{proof}[Proof of Theorem~\ref{main2}]
Note that if we replace $f(X)$ by $\rho(X)\circ f(X)\circ\eta(X)$ for any degree-one $\rho(X),\eta(X)\in\F_q[X]$, then we do not change the degree of $f(X)$, the permutation property of $f(X)$, or the number of elements $a\in\F_q$ for which $g_a(X)$ permutes $\F_q$.
%
%
Thus we may make such a replacement without affecting the truth of the hypotheses of Theorem~\ref{main2}, or of each of (1)--(3). Hence if (3) holds then Lemma~\ref{easy} implies that (2) holds. Plainly (2) implies (1), 
so it remains to show that (1) implies (3).

Henceforth we assume that (1) holds. Upon replacing $f(X)$ by $\rho(f(X))$ for some degree-one $\rho(X)\in\F_q[X]$, we may assume that $f(X)$ is monic and $f(0)=0$. Let $\Lambda$ be the set of all $a\in\F_q$ for which $g_a(X)$ permutes $\F_q$. By Lemma~\ref{PPEP}, both $f(X)$ and each $g_a(X)$ with $a\in \Lambda$ are exceptional polynomials.

We first show that $f(X)$ is not separable. For each $a\in \Lambda$, the polynomial $h_a(X):=Xg_a(X)$ equals $f(X+a)+f(a)$, and hence permutes $\F_q$. Thus the fact that $h_a(0)=0$ implies that $g_a(X)$ has no nonzero roots in $\F_q$. Since $g_a(X)$ permutes $\F_q$, it follows that $g_a(0)=0$. Thus $0$ is a multiple root of $h_a(X)$, so that $0=h_a'(0)=f'(a)$.
Since $f(X)$ permutes $\F_q$, the values $f(a)$ with $a\in \Lambda$ are pairwise distinct, so that the number of critical values of $f(X)$ is at least $\abs{\Lambda}\ge N$, and hence is greater than $\log_3(n)$. Now Corollary~\ref{EPbp} implies that $f(X)$ is not separable.

Write $f(X):=\widetilde{f}(X^Q)$ where $\widetilde{f}(X)\in\F_q[X]\setminus\F_q[X^2]$ and $Q:=2^{\ell}$ with $\ell\ge 1$. Then $\widetilde{f}(X)$ is monic and $\widetilde{f}(0)=0$. If $\widetilde{f}(X)=X$ then by Lemma~\ref{easy} we obtain (3) with $n=2^{\ell}$. 
Henceforth we assume that $m:=\deg(\widetilde{f})$ is at least $2$. Since $\widetilde{f}(X)$ permutes $\F_q$, Lemma~\ref{PPEP} implies that $\widetilde{f}(X)$ is exceptional. Since $\widetilde{f}(X)$ is separable, we have $m\ne 2$ by Lemma~\ref{34}, so that $m\ge 3$ and thus $n=mQ\ge 6$. 
Since $\widetilde{f}(X)$ is separable and exceptional of degree $m$, Corollary~\ref{EPbp} implies that $\widetilde{f}(X)$ has at most $\log_3(m)$ critical values. 
Since $f(X)$ permutes $\F_q$, it follows that there are at most $\log_3(m)$ elements $a\in\F_q$ for which $f(a)$ is a critical value of $\widetilde{f}(X)$. Let $\Gamma$ be the set of elements $a\in \Lambda$ for which $f(a)$ is not a critical value of $\widetilde{f}(X)$, so that 
\[
\abs{\Gamma}\ge \abs{\Lambda}-\log_3(m) > 3.
\]

If $a\in \Lambda$ then $g_a(X)\in X\F_q[X^2]$, so by Corollary~\ref{oddEP} we have $g_a(X)=v_a(X)\circ u_a(X)$ where $u_a(X)\in \mathcal{M}$ and $v_a(X)$ is a composition of polynomials in $\mathcal{M}$, where $\mathcal{M}$ consists of all $X^s$ for odd primes $s$ and all $D_s(X,b)$ with $b\in\F_q^*$ and prime $s\ge 5$. For each $a\in \Lambda$ we fix one such choice of $u_a(X)$ and $v_a(X)$. 
By Lemma~\ref{Dicksonram}, if $u(X)\in \mathcal{M}$ then $0$ is the unique critical value of $u(X)$. Moreover, the ramification multiset $\R_u(0)$ is either $[\deg(u)]$ or $[1,2^*]$, where in either case $0$ is the unique root of $u(X)$ with odd multiplicity. By Lemma~\ref{ram}, it follows that any composition of finitely many polynomials in $\mathcal{M}$ has a unique root with odd multiplicity, and has no roots with multiplicity divisible by $4$. 

For any $a\in \Gamma$, the polynomial $\widetilde{f}(X)+f(a)$ is squarefree, so that all roots of $\widetilde{f}(X^Q)+f(a)$ in $\mybar\F_q$ have multiplicity $Q$, whence $\R_{g_a}(0)=[Q-1,Q^*]$. Since $g_a(X)$ is a composition of polynomials in $\mathcal{M}$, we know that no elements of $\R_{g_a}(0)$ are divisible by $4$, so we must have $Q=2$. 
Thus $g_a(X)$ is a composition of polynomials of the form $D_s(X,b)$ with $b\in\F_q^*$ and $s\ge 5$ prime. Since $\deg(g_a)=n-1$, it follows that $n\not\equiv 1\pmod 3$. Then Lemma~\ref{elem} implies that $f(X)=\rho(X)\circ X^6\circ\eta(X)$ for some degree-one $\rho(X),\eta(X)\in\F_q[X]$, so Lemma~\ref{easy} implies that $k$ is odd and thus item (3) of Theorem~\ref{main2} holds.
\end{proof}

\begin{lemma}\label{elem}
Let $q=2^k$ with $k\ge 1$, and pick a non-separable $f(X)\in\F_q[X]$ of degree $n\ge 6$. Let $\Delta$ be the set of all $a\in\F_q$ for which $g_a(X):=(f(X+a)+f(a))/X$ equals $v(X)\circ D_{s}(X,b)$ for some $v(X)\in\F_q[X]$, some integer $s\ge 5$, and some $b\in\F_q^*$. If $\abs{\Delta}>3$ then $f(X)=\rho(X)\circ X^6\circ\eta(X)$ for some degree-one $\rho(X),\eta(X)\in\F_q[X]$.
\end{lemma}

\begin{proof}
The hypotheses and conclusion are unchanged if we replace $f(X)$ by $\rho(f(X))$ for some degree-one $\rho(X)\in\F_q[X]$, so we make such a replacement in order to assume that $f(X)$ is monic and $f(0)=0$.
For each $a\in \Delta$, make one choice of $v_a(X)$, $s_a$, and $b_a$ satisfying the conditions in Lemma~\ref{elem}. For ease of notation, we write $s$ for $s_a$ and $r$ for $\deg(v_a)$, although both $s$ and $r$ might depend on $a$. Then $rs=n-1$ is odd, so both $r$ and $s$ are odd; thus
\[
4\mid (r+1)(s+1)=n+r+s,
\]
so that $r+s\equiv -n\equiv n\pmod 4$. Since $f(X)$ is non-separable, we have $f(X)\in\F_q[X^2]$. Thus all terms of $g_a(X)$ have odd degree, so since the same is true of $D_s(X,b_a)$, it follows that the same is also true of $v_a(X)$.
Write $f(X)=X^n + c X^{n-2} + d X^{n-4} + O(X^{n-6})$, where $O(X^t)$ denotes a polynomial in $\F_q[X]$ of degree at most $t$ (with the convention that the zero polynomial has degree $-\infty$). Then
\[
g_a(X) = X^{n-1} + \Bigl(c+\frac{n}2 a^2\Bigr) X^{n-3} + \Bigl(d + \frac{n-2}2 c a^2 + \frac{n(n-2)}8 a^4\Bigr) X^{n-5} + O(X^{n-7}).
\]
For each $a\in \Delta$, by \eqref{dickson} we have
\[
D_s(X,b_a) = X^s + b_a X^{s-2} + \frac{s-3}2 b_a^2 X^{s-4} + O(X^{s-6}),
\]
and thus
\[
D_s(X,b_a)^r = X^{rs} + b_a X^{rs-2} + \Bigl(\frac{r-1}2 b_a^2+\frac{s-3}2 b_a^2\Bigr) X^{rs-4} + O(X^{rs-6}).
\]
Since
\[
g_a(X) = v_a(X) \circ D_s(X,b_a) = D_s(X,b_a)^r + O(X^{rs-2s}),
\]
we may equate degree-$t$ terms in $g_a(X)$ and $D_s(X,b_a)^r$ for each $t\in\{rs-2,rs-4\}$ to conclude that
\begin{equation}\label{cval}
c+\frac{n}2 a^2 = b_a
\end{equation}
and
\begin{equation}\label{dval}
d + \frac{n-2}2 c a^2 + \frac{n(n-2)}8 a^4 = \frac{r+s}2 b_a^2=\frac{n}2 b_a^2.
\end{equation}
Substituting the value of $b_a$ from \eqref{cval} into \eqref{dval} yields
\[
d + \frac{n-2}2 c a^2 + \frac{n(n-2)}8 a^4 = 
\frac{n}2\Bigl(c^2+\frac{n}2 a^4\Bigr),
\]
which may be rearranged to yield
\[
\Bigl(d+\frac{n}2 c^2\Bigr) + \frac{n-2}2 c a^2 + \Bigl(\frac{n(n-2)}8+\frac{n^2}4\Bigr) a^4 = 0.
\]
Thus $a^2$ is a root of $G(X)$ where
\[
G(X):=\Bigl(d+\frac{n}2 c^2\Bigr) + \frac{n-2}2 c X + \Bigl(\frac{n(n-2)}8+\frac{n^2}4\Bigr) X^2.
\]
Since $\abs{\Delta}>2$, and any two elements of $\Delta$ have distinct squares since $q$ is even, we see that $G(X)$ has more than two roots. But $\deg(G)\le 2$, so that $G(X)=0$. Setting the coefficients of $G(X)$ to $0$ yields
\begin{align*}
\frac{n}2 c^2&=d,\\
\frac{n-2}2 c &= 0,
\intertext{and}
\frac{n(n-2)}8 &\equiv \frac{n^2}4\pmod 2.
\end{align*}
Writing $m:=n/2$, the third equation says that $m^2+m\equiv 0\pmod 4$, so that $m$ is congruent to $0$ or $3$ mod $4$, whence $n$ is congruent to $0$ or $6$ mod $8$. If $n\equiv 0\pmod 8$ then $c=0$ and thus $b_a=c+\frac{n}{2}a^2=0$, contradiction. Thus $n\equiv 6\pmod 8$, so that $d=c^2$.

Henceforth we will assume $n\ne 6$, since if $n=6$ then we have shown that $f(X)=X^6+cX^4+c^2 X^2=(X+\sqrt{c})^6+c^3$, so the conclusion of Lemma~\ref{elem} holds. Write $f(X)=X^n + cX^{n-2}+c^2 X^{n-4}+eX^{n-6}+O(X^{n-8})$. Recalling that $n\equiv 6\pmod 8$, for $a\in \Delta$ we compute that the coefficient of $X^{n-7}$ in $g_a(X)$ is
\begin{equation}\label{co7}
a^6 + c a^4 + c^2 a^2 + e.
\end{equation}
We have
\[
D_s(X,b_a) = X^s + b_a X^{s-2} + \frac{s-3}2 b_a^2 X^{s-4} + \frac{s-5}2 b_a^3 X^{s-6} + O(X^{s-8})
\]
(note that this is valid when $s=5$, since then the coefficient of $X^{s-6}$ is $0$).
The coefficient of $X^{rs-6}$ in $D_s(X,b_a)^r$ is
\begin{equation}\label{coeff}
\frac{r-1}2 b_a^3 + \frac{s-5}2 b_a^3.
\end{equation}
Since $(v_a(X)-X^r)\circ D_s(X,b_a)$ has degree at most $(r-2)s\le rs-10$, the coefficient of $X^{rs-6}$ in $g_a(X)=v_a(X)\circ D_s(X,b_a)$ is given by \eqref{coeff}. Thus we may equate \eqref{co7} and \eqref{coeff} to obtain
\begin{equation}\label{equ}
a^6 + ca^4 + c^2 a^2 + e = \frac{r+s-6}2 b_a^3.
\end{equation}
Since $r+s\equiv n\equiv 2\pmod 4$, the right side of \eqref{equ} vanishes.  Thus $a^2$ is a root of
\[
H(X):=X^3 + cX^2 + c^2 X + e.
\]
Since $\abs{\Delta}>3$, and squaring is injective on $\Delta$, it follows that $H(X)$ has more than $3=\deg(H)$ roots, which is impossible.
\end{proof}


\section{Elementary proof} \label{sec:elem}

Since our proof of Theorem~\ref{main} used consequences of the classification of finite simple groups, it may be of interest to prove a variant of Theorem~\ref{main} using elementary methods.
In this section we sketch an elementary proof of the variant of Theorem~\ref{main} in which $N$ is replaced by $n/2$. The proof closely follows the proof of Theorem~\ref{main}, with appropriate modifications at certain spots.

As in the previous section, it suffices to prove the variant of Theorem~\ref{main2} in which $N$ is replaced by $n/2$. Moreover, as before, item (3) of the variant of Theorem~\ref{main2} implies (2), and (2) implies (1), so we will assume (1) and deduce (3).
Again we may assume that $f(X)$ is monic and $f(0)=0$, and we define $\Lambda$ as before. In order to show that $f(X)$ is separable, as in the proof of Theorem~\ref{main2} we see that the number of critical values of $f(X)$ is at least $\abs{\Lambda}$, which in this case is at least $n/2$, so that $f(X)$ is not separable by Lemma~\ref{rh}. 

Defining $\widetilde{f}(X)$, $Q$, and $m$ as before, again we may assume that $m\ge 3$ and $n\ge 6$. Lemma~\ref{rh} implies that $\widetilde{f}(X)$ has at most $(m-1)/2$ critical values; since $\widetilde{f}(X)$ permutes $\F_q$, it follows that there are at most $(m-1)/2$ elements $a\in\F_q$ for which $f(a)$ is a critical value of $\widetilde{f}(X)$. Defining $\Gamma$ as before, we have
\[
\abs{\Gamma}\ge \abs{\Lambda}-\frac{m-1}2 \ge \frac{n+2}4
\]
since $\abs{\Lambda}\ge n/2\ge m$. 
Defining $u_a(X)$ and $v_a(X)$ as before, the argument in the proof of Theorem~\ref{main2} shows that $Q=2$ and $n\not\equiv 1\pmod 3$. Moreover, each $g_a(X)$ with $a\in\Gamma$ is a composition of polynomials of the form $D_s(X,b)$ with $b\in\F_q^*$ and $s\ge 5$ prime.

Next we show that $n=6$. First, Lemma~\ref{elem} implies $n\le 10$, so that $n\in\{6,8,10\}$. Second, we have $n\ne 10$ since $n\not\equiv 1\pmod 3$. Finally, if $n=8$ then $m=n/Q=4$, so that $\widetilde{f}(X)$ is a degree-$4$ separable exceptional polynomial. By Lemma~\ref{34} $\widetilde{f}(X)$ has no critical values, so that $\abs{\Gamma}=\abs{\Lambda}\ge n/2=4$, which by Lemma~\ref{elem} leads to a contradiction. 

Pick $a\in \Gamma$, we have $g_a(X) = D_5(X,b)$ for some $b\in\F_q^*$, which implies
\[
f(X+a)+f(a) = Xg_a(X) = X D_5(X,b) = (X+\sqrt{b})^6+b^3, 
\]
so that $f(X) = \rho(X)\circ X^6\circ\eta(X)$ for some degree-one $\rho(X),\eta(X)\in\F_q[X]$. Thus Lemma~\ref{easy} implies that $k$ is odd, so that item (3) of holds.


\section{Heuristic}\label{sec:heuristic}

In this section we give a heuristic which predicts that, if various unknown functions $\F_q\to\F_q$ are roughly as likely to be permutations as would a randomly chosen function $\F_q\to\F_q$, then for any large $q$ there would be $(q+2)$-element subsets $\mathcal S$ of $\bP^2(\F_q)$
which have a subset $\mathcal T$ of size roughly $\log_3 q$ such that no line contains a point of $\mathcal T$ and at least two other points of $\mathcal S$.
By Lemmas~\ref{o1} and \ref{o2}, this says that ``at random'' there would be polynomials $f(X)\in\F_q[X]$ which permute $\F_q$ and for which there are roughly $\log_3 q$ elements $a\in\F_q$ such that $g_a(X):=(f(X+a)-f(a))/X$ permutes $\F_q$.

We first show that there is no loss in restricting to polynomials $f(X)$ of degree less than $q$. Each function $\varphi\colon\F_q\to\F_q$ is induced by a unique polynomial $f(X)\in\F_q[X]$ of degree less than $q$, and moreover if $f(X)$ permutes $\F_q$ and $q>2$ then $f(X)$ cannot have degree $q-1$. 
The function $\varphi$ is also induced by infinitely many other polynomials, all of which are congruent mod $X^q-X$. Next we note that the function $g_a\colon\F_q^*\to\F_q$ is determined by $\varphi$. 
We showed in the proof of Lemma~\ref{o2} that, if $q>2$ and $\deg(f)<q-1$, then $g_a(X)$ is injective on $\F_q^*$ if and only if $g_a(X)$ permutes $\F_q$. Thus, if $f(X)\equiv f_1(X)\pmod{X^q-X}$, and $\deg(f)<q$, then $(f(X+a)-f(a))/X$ permutes $\F_q$ for any $a\in\F_q$ such that $(f_1(X+a)-f_1(a))/X$ permutes $\F_q$. Hence we may assume $\deg(f)<q$.

Now, there are $q!$ permutations of $\F_q$, each of which corresponds to a unique polynomial $f(X)$ of degree less than $q$, and we want to show that if $q$ is large then we expect some such $f(X)$ to allow at least $N\approx\log q$ out of $q$ associated functions $g_a(X)$ to permute $\F_q$.  
We will assume that each $g_a(X)$ is roughly as likely to permute $\F_q$ as is a randomly chosen function $\F_q\to\F_q$, and moreover that if $a\ne b$ then the permutation property of $g_a(X)$ is independent of the permutation property of $g_b(X)$.
Write $M$ for the expected number of pairs $(f(X),\mathcal T)$, where $f(X)$ is a polynomial in $\F_q[X]$ of degree less than $q$ which permutes $\F_q$, and $\mathcal T$ is a set of $N\approx\log_3(q)$ elements $a$ in $\F_q$, where $g_a(X)$ permutes $\F_q$ for each $a\in\mathcal T$. Then $M$ is at least roughly
\begin{equation}\label{E}
q!\binom{q}{N} \Bigl(\frac{q!}{q^q}\Bigr)^N.
\end{equation}
Using Stirling's approximation $q!\approx \sqrt{2\pi q}\,(q/e)^q$, and approximating $\binom{q}{N}$ by $q^N/N!\approx q^N (e/N)^N/\sqrt{2\pi N}$, we see that \eqref{E} is approximately
\begin{multline*}
\sqrt{2\pi q}\,\Bigl(\frac{q}e\Bigr)^q\cdot \Bigl(\frac{qe}N\Bigr)^N (2\pi N)^{-1/2}\cdot \Bigl(\frac{\sqrt{2\pi q}}{e^q}\Bigr)^N = \frac{q^{q+(3N+1)/2}\cdot e^{-q+N-Nq}\cdot\sqrt{2\pi}^N}{N^{N+1/2}} \\
= e^{-q+N-Nq+(q+(3N+1)/2)\log q - (N+1/2)\log N + N\log(\sqrt{2\pi})}.
\end{multline*}
Since $N\approx \log_3 q$, which equals $(\log q)/c$ where $c:=\log 3$ is greater than $1$, the right side of the above equality is approximately
\[
e^{-q+N+Nq(c-1)
+cN(3N+1)/2 - (N+1/2)\log N + N\log(\sqrt{2\pi})},
\]
which becomes large when $q$ is large.

\begin{rmk}
This heuristic applies regardless of whether $q$ is even or odd. However, Proposition~\ref{oddq} shows that the conclusion of the heuristic is not true when $q$ is odd. 
The reason for this will become clear in the proof of Proposition~\ref{oddq}, which shows that there is a causal reason why $g_a(X)$ cannot permute $\F_q$, thus violating the heuristic's assumption that $g_a(X)$ has probability $q!/q^q$ of permuting $\F_q$.
However, when $q$ is even we know no reason why the heuristic's assumption would be violated for most choices of $f(X)$.
\end{rmk}


\section{Proof of Proposition~\ref{oddq}}\label{sec:odd}

In this section we prove Proposition~\ref{oddq}. In light of Lemmas~\ref{o1} and \ref{o2}, it suffices to show the following. 

\begin{prop}\label{oddq2}
If $q$ is an odd prime power, and $f(X)\in\F_q[X]$ permutes\/ $\F_q$, then there is no element $a\in\F_q$ for which $\displaystyle{g_a(X):=\frac{f(X+a)-f(a)}X}$ permutes\/ $\F_q$.
\end{prop}

\begin{proof}
Suppose to the contrary that there is some $a\in\F_q$ for which both $f(X)$ and $g_a(X)$ permute $\F_q$. Then $h_a(X):=X g_a(X)$ equals $f(X+a)-f(a)$, and hence permutes $\F_q$. 
Since $h_a(X)$ permutes $\F_q$ and fixes $0$, we see that $h_a(X)$ has no roots in $\F_q^*$, so that also $g_a(X)$ has no roots in $\F_q^*$, and since $g_a(X)$ permutes $\F_q$ we conclude that $g_a(0)=0$. Since $g_a(X)$ and $h_a(X)$ permute $\F_q$ and fix $0$, they each permute $\F_q^*$. Now we compute
\[
\theta:=\prod_{b\in\F_q^*} h_a(b)
\]
in two ways. Recall that $\prod_{b\in\F_q^*} b=-1$. Thus, since $g_a(X)$ and $h_a(X)$ permute $\F_q^*$, we have
\[
\theta=\prod_{b\in\F_q^*} b = -1
\]
and
\[
\theta=\prod_{b\in\F_q^*} \Bigl(b g_a(b)\Bigr) = \prod_{b\in\F_q^*} b \cdot\prod_{c\in\F_q^*} g_a(c) = \prod_{b\in\F_q^*} b \cdot \prod_{c\in\F_q^*} c = 1,
\]
whence $1=-1$ in $\F_q$, contradiction.
\end{proof}

\begin{rmk}
The only place in the above proof where we use the hypothesis that $q$ is odd is in order to deduce a contradiction in the last step from the equality $1=-1$.
\end{rmk}



\end{document}